\documentclass[11pt]{amsart}
\usepackage{latexsym}
\usepackage{amsfonts}

\usepackage{hyperref}
\usepackage[top=1.1in, bottom=1in, left=1in, right=1in]{geometry}

\usepackage{amsmath,amsthm,amssymb}

\newtheorem{definition}{Definition}

\newtheorem{proposition}{Proposition}[section]
\newtheorem{theorem}[proposition]{Theorem}

\newtheorem{corollary}{Corollary}

\newcommand{\Z}{\mathbb{Z}}

\newcommand{\Q}{\mathbb{Q}}

\newcommand{\F}{\mathbb{F}}

\begin{document}

\author[]{Vladimir Shpilrain}
\address{Department of Mathematics, The City College of New York, New York,
NY 10031} \email{shpilrain@yahoo.com}

\title[Complexity of algorithmic problems]{Complexity of some algorithmic problems in groups:\\
a survey}

\begin{abstract}
In this survey, we address the worst-case, average-case, and generic-case time complexity of the word  problem and some other algorithmic problems in several classes of groups and show that it is often the case that the average-case complexity of the word problem is linear with respect to the length of an input word, which is as good as it gets if one considers groups given by generators and defining relations. At the same time, there are other natural algorithmic problems, for instance, the geodesic (decision) problem or Whitehead's automorphism problem, where the average-case time complexity can be sublinear, even constant.

\end{abstract}

\maketitle

%\hfill {\small To Gerhard Rosenberger on his 80th birthday}

\section{Introduction}\label{intro}

Algorithmic problems for infinite groups have received a lot of attention over the last 100+ years, since Dehn's work back in 1911. More recently, with the increasing interest in computational complexity, many algorithmic problems of group theory have been reexamined with regard to efficiency.

The worst-case complexity of group-theoretic algorithms has been studied probably since small cancellation groups were introduced and it was noticed that the word problem in these groups admits an efficient solution by Dehn's algorithm, see e.g. \cite{LS}.
Time complexity of Dehn's algorithm was eventually shown to be linear (with respect to the ``length" of the input).
The interest in complexity of group-theoretic algorithms was afterwards gradually increasing over the rest of the 20th century, and then experienced explosive growth in the 21st century, in part due to input from theoretical computer science (notably the Knuth-Morris-Pratt algorithm and the Master theorem for divide-and-conquer recurrences), and in part due to the introduction of the generic-case and the average-case complexity into group theory about 20 years ago.

Genericity of group-theoretic properties and generic-case complexity of group-theoretic algorithms were introduced in \cite{AO} and \cite{KMSS}. Average-case complexity, formally introduced in \cite{Levin} (although informally considered before in \cite{Knuth}), was addressed in the context of group theory for the first time in \cite{KMSS2}. Specifically, the authors of \cite{KMSS2} addressed the average-case complexity of the word and subgroup membership problems in non-amenable groups and showed that this complexity was often linear. The fact that the considered groups were non-amenable played an essential role in computing the average-case complexity due to the fact that the cogrowth rate of any non-amenable group is strictly smaller than the growth rate of the ambient free group, see \cite{Cohen}, \cite{Grigorchuk}. %We give more details on the cogrowth function of a group in Section \ref{Average-casesolvable}.
In contrast, most groups considered in a recent paper \cite{OSh} are amenable. The case of amenable groups is more difficult; general tricks of \cite{KMSS2} are not applicable in this case, and one has to actually investigate the structure of groups in question.

There are several ways to define the average-case time complexity of a group-theoretic algorithm $\mathcal{A}$ that takes words as inputs; below is one of the natural ways.
Let $W_n$ denote the set of all words of length $n$ in a finite group alphabet. For a word $w \in W_n$, let $T(w)$ denote the time that the algorithm $\mathcal{A}$ (for a classical Turing machine) works on input $w$.

One can then define the average-case time complexity of the algorithm $\mathcal{A}$ on inputs of length $n$ as
\vskip -0.5cm

\begin{equation}\label{avcase1}
\frac{1}{|W_n|} \sum_{w\in W_n} T(w).
\end{equation}

\subsection{Why the average-case complexity is often strictly lower than the worst-case complexity} \label{Vegas}

On the intuitive level, the reason is that inputs $w\in W_n$ for which $T(w)$ is high are ``sparse".
If one is able to quantify this ``sparsity", then one can try to split $W_n$ in a disjoint union $W_n=\cup_j V_n^{(j)}$ of sets $V_n^{(j)}$, so that the formula (\ref{avcase1}) is stratified as
\vskip -0.5cm

\begin{equation}\label{avcase2}
\frac{1}{|W_n|} \sum_j |V_n^{(j)}| \cdot T(w \in V_n^{(j)}) = \sum_j \frac{|V_n^{(j)}|}{|W_n|} \cdot T(w \in V_n^{(j)}).
\end{equation}

This formula is not very practical though because the number of summands may be too large.
Keeping in mind that our goal is typically reduced to finding an {\it upper bound} for the average-case complexity rather  than its precise value, we will use a slightly different stratification that is more practical.
Denote by $\overline{T}(w \in V_n^{(j)})$ an upper bound on $T(w)$ for $w \in V_n^{(j)}$ and use the following sum instead of (\ref{avcase2}):
\vskip -0.5cm

\begin{equation}\label{avcase3}
\sum_j \frac{|V_n^{(j)}|}{|W_n|} \cdot \overline{T}(w \in V_n^{(j)}).
\end{equation}

In this sum, the number of summands can be rather small (typically 2 or 3, although it can be larger in some cases),
depending on an upper bound on the average-case complexity one would like to establish.

Thus, our strategy will be to find suitable upper bounds on $\frac{|V_n^{(j)}|}{|W_n|}$ for those $V_n^{(j)}$ where $\overline{T}(w \in V_n^{(j)})$ is high.

An alternative strategy (used in \cite{KMSS2} and in \cite{White}, say) is, for a given input, to run two algorithms in parallel. One algorithm, call it {\it honest}, always terminates in finite time and gives a correct result. The other algorithm, a {\it Las Vegas algorithm}, is a fast randomized algorithm that never gives an incorrect result; that is, it either  produces the correct result or informs about the failure to obtain any result. (In contrast, a {\it Monte Carlo algorithm} is a randomized algorithm whose output may be incorrect with some (typically small) probability.)

A Las Vegas algorithm can improve the time complexity of an honest, ``hard-working", algorithm that always gives a correct answer but is slow. Specifically, by running a fast
Las Vegas algorithm and a slow honest algorithm in parallel, one often gets an algorithm that always terminates with a correct answer and whose average-case complexity is somewhere in between. This idea was used in \cite{KMSS2} where it was shown, in particular, that if a group $G$ has the word problem solvable in subexponential time and if $G$ has a non-amenable factor group where the word problem is solvable in a complexity class $\mathcal{C}$, then there is an honest algorithm that solves the word problem in $G$ with average-case complexity in $\mathcal{C}$.

\section{Background}\label{previous}

Bibliography on the worst-case complexity of various group-theoretic problems (the word problem, conjugacy problem, subgroup membership problem, geodesic problem, equations in groups) includes hundreds of papers, so a survey on those would take a hefty monograph. Perhaps just as a few milestones, we can mention here papers \cite{Birget}, \cite{BKL}, \cite{Vershik}, \cite{PatRaz},  \cite{Saul} and the monograph  \cite{Lohrey}.

The shift of the interest, in the last 20 years or so, from the worst-case to the generic-case and average-case complexity is in part due to applications of non-abelian groups in information security, see \cite{book2}, \cite{book3}, \cite{RamaNotices}, \cite{bookDK} and, more generally, due to emerging connections between group theory and theoretical computer science (see e.g. \cite{bookof10} for a general overview).

Some ideas in the spirit of the average-complexity in group theory were actually put forward by Gromov in his essay \cite{Gromov}. Specifically, he proposed to consider the {\it averaged} Dehn function instead of the usual Dehn function and suggested that the former might be asymptotically smaller than the latter. Relevance of the averaged Dehn function to the average-case complexity of the word problem was pointed out by Young \cite{Young}: the averaged Dehn function represents the complexity of verifying that an element of a group represents the identity, while the technique of  \cite{KMSS2} established low average-case complexity in many (non-amenable) groups by showing that in those groups, most elements are {\it not} the identity.

We also mention here an earlier paper by Wang \cite{Wang95} with ``average-case complexity" in the title. However, Wang's approach is very different in substance: in his approach the set of instances of a problem involves all finite presentations of a group,
rather than just a fixed one.

Below we first briefly survey the approach of \cite{KMSS2} applicable to non-amenable  groups, and then focus on a more difficult case of amenable groups.

\subsection{Non-amenable groups}

The basic idea of \cite{KMSS2} is very straightforward and is often used in
practice.  If we have a total algorithm $\Omega_1$ solving a decision
problem $\mathcal{D}$ whose worst-case complexity is ``not too high''
and we also have a partial algorithm $\Omega_2$ solving the problem with ``very low''  generic-case complexity, then by running  $\Omega_1$ and  $\Omega_2$
in parallel we get a total algorithm  $\Omega_1 ||  \Omega_2$
for which we can often prove to have low average-case complexity.

To formalize this idea, let us start with a definition. Below $X^*$ denotes the set of all group words in the alphabet $X$.

\begin{definition}
Let $X$ be a finite alphabet with $k\ge 2$ elements. A
\emph{discrete probability measure} on $X^*$ is a function $\mu:
X^*\to [0,1]$ such that $\sum_{w\in X^*} \mu(w)=1$.

We  say that $\mu$ is \emph{length-invariant} if for any
words $w,w'\in X^*$ with $|w|=|w'|$ we have $\mu(w)=\mu(w')$.

We say that a measure $\mu: X^*\to [0,1]$ is \emph{tame} if there
exists a subexponential
function $g(n)\ge 0$ such that for any $w\in X^*$ we have $\mu(w)\le
\frac{g(|w|)}{k^{|w|}}$.
\end{definition}

Then the main result of \cite{KMSS2} is

\begin{theorem}\label{AB}
Let $G$ be a finitely presented group where the word problem
is in $\mathcal SubExp$ (i.e., is solvable in subexponential time).
Suppose $G$ has a subgroup of finite index which has a
non-amenable quotient group $\overline  G_1$ whose  word problem is
solvable in a complexity class $\mathcal C$, where
 $\mathcal{C} \subseteq \mathcal SubExp$.

Then the word problem in $G$ is solvable with average-case complexity in
$\mathcal C$, uniformly relative to the family of all
length-invariant discrete probability measures $\mu: X^*\to [0,1]$.

Moreover, for any tame discrete probability measure $\mu$ on $X^*$ the word problem in $G$ is solvable with average-case complexity in
$\mathcal C$ relative to $\mu$.

\end{theorem}

Theorem~\ref{AB} applies to a surprisingly wide range of examples, although the idea behind the proof of this theorem is quite simple.
We run in parallel the total subexponential algorithm for the word problem of
$G$ and the partial generic algorithm provided by the quotient group
$\overline G_1$.  If  the image of an element of $g$ is nontrivial in
$\overline G_1$, then the element $g$ is obviously nontrivial in $G$. On all but an
exponentially negligible set of inputs this  partial  algorithm will
actually terminate within the complexity bound $\mathcal C$ and it
turns out that the average-case complexity of the combined
algorithm is in $\mathcal C$.

Below are some more concrete corollaries of Theorem~\ref{AB}.

\begin{corollary}\label{BB}
Let $G$ be a finitely presented group where the word problem is
solvable in subexponential time. Let $A$ be a generating set of
$G$ and let $X=A\cup A^{-1}$. Let $\Re$ be the family of all
length-invariant discrete probability measures $\mu:X^*\to [0,1]$.

\begin{enumerate}
\item Suppose that $G$ has a subgroup of finite index that
has a non-elementary word-hyperbolic
quotient group. Then the word problem in $G$ is
solvable with linear time average-case complexity uniformly
relative to $\Re$.

Moreover, for any tame discrete probability measure $\mu$ on $X^*$ the
word problem in $G$ is solvable with linear time average-case complexity
 relative to $\mu$.

\item Suppose that $G$ has a subgroup of finite index that has
a non-amenable automatic quotient group.
Then the word problem in $G$ is
solvable with quadratic time average-case complexity uniformly
relative to $\Re$.

Moreover, for any tame discrete probability measure $\mu$ on $X^*$ the
word problem in $G$ is solvable with quadratic time average-case
complexity relative to $\mu$.
\end{enumerate}
\end{corollary}

This corollary applies, in particular, to braid groups $B_n$ since $B_n$ has a subgroup of finite index, the pure braid group $P_n$, that admits a homomorphism onto the free group $F_2$. Since the word problem for $B_n$ is solvable in quadratic time (see e.g. \cite{BKL}), Corollary \ref{BB} implies that the word problem in $B_n$ has linear time average-case complexity.

Another application is to \emph{Artin groups of extra large
type}.  An Artin group $G$ has a presentation
\begin{equation*}\label{art}
G=\langle a_1,\dots,  a_t \, | u_{ij}=u_{ji}, \text{ for }1\le
i<j\le t \rangle
\end{equation*}
where for $i\ne j$
\[
u_{ij}:=\underbrace{a_ia_ja_i\dots}_{m_{ij} \text{ terms }}
\]
and where $m_{ij}=m_{ji}$ for each $i<j$. We allow $m_{ij}=\infty$
in which case the relation $u_{ij}=u_{ji}$ is omitted from the
above presentation. $G$ is of \emph{extra large
type} if $t \ge 3$ and  all $m_{ij}\ge 4$.  Peifer~\cite{Pe} proved that all such groups are automatic and therefore have the word problem solvable in quadratic time. An Artin group $G$ has an associated Coxeter group $G_1$ which is the
quotient group obtained by setting the squares of the generators
$a_i$ equal to the identity. If $G$ is of extra large type then so
is its Coxeter quotient

\[
G_1=G/ncl_G(a_1^2, \dots, a_t^2).
\]

This $G_1$ is known to be a nonelementary word-hyperbolic
group if $t>2$. Thus,  by Corollary~\ref{BB}, the word problem in any  Artin group
$G$ of extra large type is solvable in linear time on average.

\subsection{Amenable groups}

For amenable groups, the trick used in the proof of Theorem~\ref{AB} may not work, as we have explained in the Introduction. Below we mention several results obtained by using different approaches.

In \cite{OSh}, it was shown that in several classes of (amenable) groups, as well as in a large class of groups that are free products, the average-case complexity of the word problem is linear. More specifically, the following results were established:

\begin{itemize}

%\smallskip
\item[1.] The average-case time complexity of the word problem in groups of matrices over integers (or rationals) is often linear; in particular, this is the case in all polycyclic groups. These results apply, in particular, to finitely generated nilpotent groups. It is also noteworthy that the worst-case complexity of the word problem in these groups is $O(n \cdot \log^2 n)$, where $n$ is the length of the input (as a word in the generating matrices). The proof uses a recent result \cite{Harvey} on fast multiplication of integers.

    For finitely generated nilpotent groups, the worst-case time complexity of the word problem is, in fact, $O(n \cdot \log^{(k)} n)$ for any integer $k\ge 1$,  where $\log^{(k)} n$ denotes the function $\log \ldots \log n$, with $k$ logarithms.
\smallskip

\item[2.] If a finitely generated group $G$ has polynomial-time worst-case complexity
of the word problem and maps onto
a (generalized) lamplighter group (the wreath product of $\Z_p$ and $\Z$), then the average-case time complexity of the word problem in $G$ is linear. This applies, in particular, to free solvable groups.
\smallskip

\item[3.] In popular groups like solvable Baumslag-Solitar groups $BS(1, n)$ and
Thompson's group $F$, the average-case time complexity of the word problem is linear.

\noindent We recall that the group $BS(1, n)$ has a presentation $BS(1, n) = \langle a, b \, | \, bab^{-1} = a^n \rangle.$ We also note, in passing, that it is a major open problem whether or not Thompson's group $F$ is amenable, see \cite{Guba} for a survey of work done on this problem.

\smallskip

\item[4.] If a finitely generated group $G$ is a free product of nontrivial groups
$A$ and $B$, both having polynomial-time worst-case complexity of the word problem, then $G$ has linear average-case time complexity of the word problem.
\smallskip

\item[5.] If a finitely generated group $G$ is a free product of nontrivial groups
$A$ and $B$, both having polynomial-time worst-case time complexity of the subgroup membership problem, then $G$ has linear average-case time complexity of the subgroup membership problem. All subgroups in question should be finitely generated.

\end{itemize}

We also mention that the problem in some sense dual to finding the average-case complexity of the word problem in a given group $G$ is finding ``non-cooperative" words $w$ for which deciding whether or not $w=1$ in $G$ is as computationally hard as it gets; these words are ``responsible" for a high worst-case complexity of the word problem in $G$. This is relevant to properties of the {\it Dehn function} of $G$, see e.g. \cite{Birget},  \cite{MU}, \cite{Young} and references therein.

The word and subgroup membership problems are not the only group-theoretic problems whose average-case complexity can be significantly lower than the worst-case complexity.
In \cite{White}, we showed that the average-case complexity of the problem of detecting a primitive element in a free group has {\it constant} time complexity (with respect to the length of the input) if the input is a cyclically reduced word. The same idea was later used in \cite{Roy} to design an algorithm, with constant average-case complexity, for detecting {\it relatively primitive} elements, i.e., elements that are primitive in a given subgroup of a free group.

%In this project, we are going to explore other group-theoretic problems where the average-case complexity has a good chance to be sublinear, even  constant, see Section \ref{goals}.

\section{Directions of research}\label{goals}

In this section, we mention several interesting (in our opinion) problems relevant to  the average-case time complexity of the word problem  in several classes of groups. We emphasize that it is often the case that the average-case complexity of the word problem is linear with respect to the length of an input word, which is as good as it gets if one considers groups given by generators and defining relations, see e.g. \cite{Sublin}.

With some other problems, e.g. the {\it geodesic (decision) problem}, it can be the case for some groups that the average-case complexity is sublinear, even  constant. %This avenue will be explored, too.

More specifically, we suggest the following directions of research. Some of them come from the ``Open problems" section of \cite{OSh}.

\medskip

%\smallskip
\noindent {\bf Direction 1.} Prove that the average-case time complexity of the word problem
in any finitely generated metabelian group is linear.

We note that in \cite{KMS}, the authors showed that there are finitely generated solvable groups of derived length 3 (even residually finite ones) with super-exponential worst-case complexity of the word problem. The average-case time complexity of the word problem in such groups cannot be polynomial since the super-exponential runtime, even on a very small set of words, will dominate everything else.

With metabelian groups, the situation is different since it is known that, unlike in solvable groups of derived length $\ge 3$ \cite{OK}, the word problem in any finitely presented metabelian group is solvable. This follows from the fact that these groups are residually finite, but the first explicit algorithm was offered in \cite{Timoshenko}.

The question of the worst-case complexity of the word problem in metabelian groups is a problem of independent interest that would have to addressed along the way.
Here we just note that the worst-case complexity of the word problem in a finitely generated {\it free} metabelian  groups was shown to be $O(n \log n)$, i.e., it is quasilinear in the length $n$ of the input \cite{Vershik}.
\smallskip

\noindent {\bf Direction 2.} Prove that the average-case time complexity of the word problem
in any finitely generated group of matrices over $\Q$ is linear.

This is a natural direction prompted by the results of \cite{OSh} about complexity of the word problem in matrix groups mentioned in Section \ref{previous} of this survey.

We also note that attacking this problem will probably involve estimating the generic-case complexity of matrix multiplication. Since the latter is a problem of independent interest, we single it out as Direction 5 below.
\smallskip

\noindent {\bf Direction 3.}  Let $G$ be a finitely generated (or even a finitely presented) group where the average-case time
complexity of the word problem is linear. We plan to address the following question:  is it necessarily the case that the generic-case time complexity (i.e., complexity on ``most" inputs) of the word problem in $G$ is linear, too?

We note that the  converse is not true, see the comment to the Direction 1 above. We also note that if the sums in the formula (\ref{avcase1}) are bounded by a linear function of $n$, then for any superlinear function $f(n)$ the densities of subsets of $G$ where $T(w)>f(n)$ (for $w\in W_n$) are approaching 0 when $n$ goes to infinity. This means that if the average-case time complexity of the word problem in $G$ is linear, then the generic-case time complexity of the word problem in $G$ can be bounded by any superlinear function.

A similar argument shows that if the average-case time complexity of the word problem in $G$ is polynomial, then the generic-case time complexity of the word problem in $G$ is polynomial, too, although the degree of a polynomial may be a little higher. More accurately, if the average-case complexity is $O(n^d)$, then the generic-case complexity is $O(n^{d+\varepsilon})$ for any $\varepsilon > 0$.

Another implication that we can mention in this context is that the exponential time average-case time complexity of the word problem in $G$ implies the exponential time (as opposed to super-exponential time) worst-case complexity of the word problem in $G$. This is because super-exponential complexity of the word problem even on just a single input will dominate everything else, so the average-case complexity would be super-exponential as well.
\smallskip

\noindent {\bf Direction 4.} This direction is somewhat informal. Prompted by a result of \cite{White} mentioned at the end of Section \ref{intro}, look for other algorithmic problems where the result of \cite{languages} on exponentially negligible proportion of words with ``forbidden subwords" can lead to the constant-time average-case complexity.

    One way to generalize the result of \cite{White} about fast (on average) detection of primitive elements of a free group $F_r$ would be finding other ``orbit-blocking" words, see \cite[Problem (F40)]{Problems}. In case that problem has a positive answer for a particular word $u$, the following restricted version  of Whitehead's automorphism problem will have constant-time average-case complexity (with respect to the length of the input word $v$): given a cyclically reduced word $v$, find out if $v=\alpha(u)$ for some $\alpha\in Aut(F_r)$.

    Another problems of that kind is the {\it geodesic (decision) problem}: given a finitely presented group $G$ and a word $w$ in the generators of $G$, find out whether or not there is a shorter word that represents the same element of $G$. ``Forbidden subwords" in this case are subwords of $w$ that are longer subwords of defining relators of $G$. By a ``longer subword" of $r$ we mean a subword whose length is more than half of the length of $r$ (in particular, any defining relator is such a ``forbidden subword"). If $w$ has such a subword, then $w$ is not a geodesic word.

    Based on this simple observation, one can show, in particular, that the average-case time complexity of the geodesic problem in any hyperbolic group is constant (with respect to the length of the input word). This is because the worst-case time complexity of this problem in any hyperbolic group is at most quadratic. In fact, if the geodesic (decision) problem in a group $G$ is solvable in subexponential time in the worst case, then the average-case complexity of this problem in $G$ is linear-time.

    One can look for other algorithmic problems where the presence of ``forbidden subwords" can quickly lead to a decision, and then use this to bound the average-case complexity.

    \smallskip

\noindent {\bf Direction 5.} This direction is relevant to Direction 2, but is also of independent interest because of connections to {\it Cayley hash functions}, see \cite{BSV}. More specifically, it is relevant to a well-known problem of estimating the {\it girth} of the Cayley graph of some special 2-generator groups of matrices over $\Z_p$; the latter problem has received a lot of attention over the last 15 years or so, see e.g. \cite{AB}, \cite{BG}, \cite{BSV},  \cite{H}.

    The problem to consider here is as follows. Let $A$ and $B$ be $2 \times 2$ matrices over $\Z$ that generate a free group. Let $w(A, B)$ be a freely reduced (semi)group word of length $n$. After evaluating $w(A, B)$ as a product of matrices, we get a $2 \times 2$ matrix, call it $W$. What is the largest (by the absolute value) possible entry of $W$, over all $w(A, B)$ of length $n$, as a function of $n$? The latter function is usually (although not always) exponential, so we are looking for an answer in the form $O(s^n)$, i.e., we are looking for the base $s$ of the exponent.

    The relevance of this problem to the problem of estimating the girth of the Cayley graph of 2-generator (semi)groups of matrices over $\Z_p$ is almost obvious. If $A$ and $B$ generate a free sub(semi)group of $SL_2(\Z)$, then there cannot be any relations of the form $u(A, B)=v(A, B)$ in $SL_2(\Z_p)$ unless at least one of the entries of the matrix $u(A, B)$ or $v(A, B)$ is at least $p$. Thus, if the largest entry in a product of $n$ matrices is of the size $O(s^n)$, then the girth of the Cayley graph of our sub(semi)group of $SL_2(\Z_p)$ is $O(\log_s p)$.

    A more difficult problem is to estimate the growth of a ``random" product of $n$ matrices, each of which is $A$ or $B$. That is, each matrix in such a product is
either $A$ or $B$, with probability $\frac{1}{2}$. This is relevant to Direction 2, specifically to the generic-case complexity of matrix multiplication and therefore to the generic-case complexity of the word problem in (two-generator) groups of matrices. We discuss an approach to this problem at the end of Section \ref{ways}.

\section{Some possible approaches}\label{ways}

In this section, we describe what we believe are promising approaches along some of the research avenues mentioned in Section \ref{goals}.

In Direction 1, one will first need to address the worst-case complexity of the  word problem in finitely generated metabelian groups. There is an explicit algorithm, due to Timoshenko \cite{Timoshenko}, for solving this problem, but the worst-case complexity of this algorithm is unknown. Perhaps it was discouraging that Timoshenko's algorithm employed Gr\"obner bases (in Laurent polynomial rings), and these are not known to be efficiently constructible. However, constructing relevant Gr\"obner bases is not part of solving the word problem for a given input word. These bases can be precomputed and then used for any input word $w$. With this in mind, it is quite likely that the worst-case time complexity of the word problem in finitely generated metabelian groups is polynomial.

If this is the case, then for the average-case complexity one can use an approach from \cite{OSh} based on the following

\begin{theorem}\cite[Theorem 15.8a]{Woess} \label{lemma Woess}
Let $G$ be a polycyclic group, and suppose that all words of length $n$ (in the given generators of $G$) are sampled with equal probability.
%with a probability measure that is finitely supported and symmetric.
Then one has the following alternative:
\medskip

\noindent {\bf (a)} $G$ has polynomial growth with degree $d$ and for a random word $w$ of length $n$, the probability of the event $w = 1$ is $O(\frac{1}{n^{d/2}})$.

\medskip

\noindent {\bf (b)} $G$ has exponential growth and  the probability of the event $w = 1$
is $O(\exp(-c \cdot n^{\frac{1}{3}}))$ for some constant $c>0$.

\end{theorem}

As an example of using this result, we show below how one can establish that in any polycyclic group, the average-case complexity of the word problem is linear. Let $G$ be an infinite polycyclic group. It has a factor group $H=G/N$ that is virtually free abelian of rank $r \ge 1$.
Thus, we first check if an input word $w$ is equal to 1 modulo $N$; this takes linear time in $n$, the length of $w$. Indeed, let $w=g_1 \cdots g_{n}$ and let $K$ be an abelian normal subgroup of $H$ such that $H/K$ is finite. Then first check if $w=1$ in $H/K$; this can be done in linear time. If $w \in K$, then rewrite $w$ in generators of $K$ using the Reidemeister-Schreier procedure. This, too, takes linear time. Indeed, if we already have a Schreier representative for a product  $g_1 \cdots g_{i-1}$, then we multiply its image in the factor group $H/K$ by the image of $g_i$ and get a representative for $g_1 \cdots g_{i}$. This step is performed in time bounded by a constant since the group $H/K$ is finite. Therefore, the whole rewriting takes linear time in $n$.

The obtained word $w'$ should have length $\le n$ in generators of $K$. Then we solve the word problem for $w'$ in the abelian group $K$, which can be done in linear time in the length of $w'$, and therefore the whole solution takes linear time in $n$.

If $w=1$ modulo $N$, then we apply an algorithm for solving the word problem in $N$ with the worst-case complexity $O(n \log^2 n)$; such an algorithm exists by a result of \cite{OSh} mentioned in Section \ref{previous}. (Recall that polycyclic groups are representable by matrices over integers.) The probability that $w=1$ modulo $N$ is  $O(\frac{1}{n^{1/2}})$ by Theorem \ref{lemma Woess}.

Therefore, the expected runtime of the combined algorithm is $O(n) + O(\frac{1}{n^{1/2}}) \cdot O(n \log^2 n)  = O(n)$.

Along the same lines, one can establish linear-time average-case complexity of the word problem in many groups $G$ that have an (infinite) abelian factor group. It is important however to have a polynomial bound on the worst-case time complexity of the word problem in $G$. (In some cases, a subexponential bound may suffice.)
In particular, this is how one can approach linear-time average-case complexity of the word problem in finitely generated metabelian groups, by first establishing a polynomial bound on the worst-case complexity.

\subsection{Constant-time average-case complexity}\label{forbidden}

One might find surprising the fact that the average-case complexity of some natural group-theoretic problems can be constant. Below we illustrate how this can happen by using a special case of the Whitehead problem for a free group. The idea of using ``forbidden subwords" in estimating the average-case complexity can be fruitful in  other situations as well.

Let $F_r$ be a free group with a free generating set $x_1, \ldots, x_r$ and let $w=w(x_1, \ldots, x_r)$.
Call an element $u \in F_r$ {\it primitive} if there is an automorphism of $F_r$ that takes $x_1$ to $u$.

We will also need the definition of the Whitehead graph of an element $w \in F_r$. The Whitehead graph $Wh(w)$ of $w$ has $2r$  vertices that correspond to $x_1, \ldots, x_r, x_1^{-1}, \ldots, x_r^{-1}$. For  each occurrence of a subword  $x_i x_j$ in the  word $w \in F_r$, there
is an  edge in $Wh(w)$ that connects the vertex $x_i$ to the vertex
$x_j^{-1}$; ~if $w$ has a subword  $x_i x_j^{-1}$, then there is an
edge connecting $x_i$ to $x_j$, etc. ~There is one more edge (the
external edge): this is the edge that connects the vertex corresponding to the last
letter of $w$ to the vertex corresponding to the inverse of the
first letter.

It was observed by Whitehead himself (see \cite{Wh} or \cite{LS}) that the Whitehead graph of any cyclically reduced primitive element $w$ of length $>2$ has either an isolated edge or a cut vertex, i.e., a vertex that, having been removed from the graph together with all incident edges, increases the number of connected components of the graph.

Now call a group word $w=w(x_1, \ldots, x_r)$ {\it primitivity-blocking} if it cannot be a subword of any cyclically reduced primitive element of $F_r$. For example, if the Whitehead graph of $w$ (without the external edge) is complete (i.e., any two vertices are connected by at least one edge), then $w$ is primitivity-blocking because in this case, if $w$ is a subword of $u$, then the Whitehead graph of $u$, too, is complete and therefore does not have a cut vertex or an isolated edge.
Examples of primitivity-blocking words are: $x_1^n x_2^n \cdots x_r^n x_1$ (for any $n \ge 2$), $[x_1, x_2][x_3, x_4]\cdots [x_{n-1}, x_{n}]x_1^{-1}$ (for an even $n$), etc. Here $[x, y]$ denotes $x^{-1}y^{-1}xy$.

The usual algorithm for deciding whether or not a given element of $F_r$ is primitive is a special case of the general Whitehead algorithm that decides, given two elements $u, v \in F_r$,  whether or not $u$ can be taken to $v$ by an automorphism of $F_r$. If $v=x_1$, the worst-case complexity of the Whitehead algorithm is at most quadratic in $\max(|u|, |v|)=|u|$. Here we are going to address the average-case complexity of the Whitehead algorithm run in parallel with a ``fast check" algorithm, when applied to recognizing primitive elements of $F_r$.

A ``fast check" algorithm  $\mathcal{T}$ to test primitivity of an input (cyclically reduced) word $u$ would be as follows. Let $n$ be the length of $u$.  The algorithm  $\mathcal{T}$ would read the initial segments of $u$ of length $k$, $k=1, 2, \ldots,$ adding one letter at a time, and build the Whitehead graph of this segment, excluding the external edge. Then the algorithm would check if this graph is complete. If it is complete, then in particular it does not have a cut vertex or an isolated edge, so the input element is not primitive.

Note that the Whitehead graph always has $2r$ vertices, so checking the property of such a graph to be complete takes constant time with respect to the length of $u$, although reading a segment of $u$ of length $k$ takes time $O(k)$. If the Whitehead graph of $u$ is complete, the algorithm returns ``$u$ is not primitive". If it is not, the algorithm just stops.

Denote the ``usual" Whitehead algorithm by $\mathcal{W}$.
Now we are going to run the algorithms $\mathcal{T}$ and  $\mathcal{W}$ in parallel; denote the composite algorithm by $\mathcal{A}$. Suppose that inputs of the algorithm $\mathcal{A}$ are cyclically reduced words that are selected uniformly at random from the set of cyclically reduced words of length $n$. Then we claim that the average-case time complexity (a.k.a. expected runtime) of the algorithm $\mathcal{A}$, working on a classical Turing machine, is $O(1)$, a constant that depends on $r$ but not on $n$.

Indeed, for complexity of the algorithm $\mathcal{T}$ we can use
a result of \cite{languages} saying that the number of (freely reduced) words of length $L$ with (any number of) forbidden subwords  grows exponentially slower than the number of all freely reduced words of length $L$.

In our situation, if the Whitehead graph of a word $w$ is not complete, that means $w$ does not have at least one $x_i^{\pm 1}x_j^{\pm 1}$ as a subword. Thus, if the Whitehead graph of any initial segment of the input word $u$ is not complete, we have at least one forbidden subword. Therefore, the probability that the Whitehead graph of the initial segment of length $k$ of the word $u$ is not complete is bounded by $s^k$ for some $s, ~0<s<1$. Thus, the average time complexity of the algorithm $\mathcal{T}$ is
\vskip -0.5cm

\begin{equation}\label{T}
\sum_{k=1}^n k\cdot s^k,
\end{equation}
\vskip -0.3cm

\noindent which is bounded by a constant.

Now suppose that the Whitehead graph of the input word $u$ of length $n$ is not complete, so that we have to employ the Whitehead algorithm $\mathcal{W}$. The probability of this to happen is bounded by $s^n$ for some $s, ~0<s<1$, as was mentioned before.

Then, the worst-case complexity of the Whitehead algorithm for detecting primitivity of the input word is known to have time complexity $O(n^2)$.

Thus, the average-case complexity of the composite algorithm $\mathcal{A}$ is

\begin{equation}\label{A}
\sum_{k=1}^n k\cdot s^k + O(n^2) \cdot O(s^n),
\end{equation}

\noindent which is bounded by a constant.

These results may be generalized as follows. Call $w \in F_r$ an
$Orb(u)$-blocking word if it cannot be a subword of any cyclically reduced $v \in Orb(u)$. Suppose there are $Orb(u)$-blocking words for a  $u \in F_r$ (cf. \cite{Problems}, Problem (F40)). Then, by using the same kind of argument as above, one would be able to show that the average-case complexity (with respect to the length of the input word $v$) of the following problem is constant: given a cyclically reduced word $v$, find out if $v=\alpha(u)$ for some $\alpha\in Aut(F_r)$.

\subsection{Complexity of a product of $n$ matrices}

Finally, we comment on Direction 5 from Section \ref{goals}. To begin, let $A = A(2) = \left(
 \begin{array}{cc} 1 & 2 \\ 0 & 1 \end{array} \right) , \hskip .2cm B = B(2) = \left(
 \begin{array}{cc} 1 & 0 \\ 2 & 1 \end{array} \right).$

The largest entry in a product of $n$ matrices each of which is $A$ or $B$ is $O((1+\sqrt{2})^n)$; this was proved in \cite{BSV}. (In fact, there is given an explicit function $f(n)$ for the largest entry.) The fastest growing largest entry in a product of $n$ matrices is when such a product is of the form $(AB)^{\frac{n}{2}}$ (assuming that $n$ is even).

To estimate the {\it average} growth of entries in a random product of $n$ matrices $A$ or $B$ (where each factor is either $A$ or $B$ with probability $\frac{1}{2}$), one can try to solve the following pair of recurrence relations for, say, the entries in the first row  of a $2\times 2$ matrix:

\noindent -- with probability $\frac{1}{2}$, $~a_n=a_{n-1}$, $~b_n=2a_{n-1} + b_{n-1}$;\\
-- with probability $\frac{1}{2}$, $~a_n=a_{n-1} + 2b_{n-1}$, $~b_n=b_{n-1}$.

One can compute the ``expectation" of $a_n$ as $\frac{1}{2}(a_{n-1} + a_{n-1} + 2b_{n-1}) = a_{n-1} + b_{n-1}.$

The same for the ``expectation" of $b_n=\frac{1}{2}(2a_{n-1} + b_{n-1} + b_{n-1}) = a_{n-1} + b_{n-1}.$

Thus, ``on average" $a_n=b_n$, hence from the above, $a_n=2a_{n-1}$, whence $a_n$ grows as $2^n$, and so does $b_n$.

The variance, however, seems to be pretty large, as suggested by computer experiments.
In particular, the largest entry in a product of 1000 random matrices ranges from $10^{273}$ to $10^{287}$. The most popular largest entry is on the order of $10^{280}$, which suggests that ``generically", the largest entry grows like $O((1.905)^n)$. %Our goal here is to find a theoretical explanation of these experimental results.

If one has in mind applications to Cayley hash functions (see Section \ref{goals}, Direction 5), then one wants the generic growth to be as slow as possible. To that end, one can consider a slightly different pair of matrices: $A = A(2) = \left(
 \begin{array}{cc} 1 & 2 \\ 0 & 1 \end{array} \right) , \hskip .2cm B = B(-2) =  \left(
 \begin{array}{cc} 1 & 0 \\ -2 & 1 \end{array} \right).$

With this pair, the largest entry in a product of $n$ matrices $A$ or $B$ is when such a product is of the form $(ABBA)^{\frac{n}{4}}$ (assuming that $n$ is a multiple of 4). The growth rate of the largest entry in the matrix $(ABBA)^{\frac{n}{4}}$ is $O((\sqrt{2+\sqrt{3}})^n)$. Note that $\sqrt{2+\sqrt{3}} \approx 1.93 ~< ~1+\sqrt{2} \approx 2.41$.

We note in passing that, interestingly, entries of $(AB)^{\frac{n}{2}}$  exhibit linear(!) growth in this case.

The average growth for the largest entry in a product of $n$ matrices $A$ or $B$ in this case is $O((\sqrt{2}^n) \approx O(1.41^n)$, and in a generic product it is about  $O(1.68^n)$. We make a disclaimer that all these growth estimates for the pair $(A(2), B(-2))$ are based just on computer experiments.

We also note that Pollicott \cite{Pollicott} studied the (maximal) Lyapunov exponent $\lambda$ for random matrix products of strictly positive matrices. In our situation, the maximal Lyapunov exponent is the natural logarithm of the  maximal base of the exponent in the growth rate of the entries in a random product of $n$ matrices $A$ or $B$. In \cite{Pollicott}, an algorithm is given for estimating $\lambda$ with any desired precision. However, our matrices $A$ and $B$ are not strictly positive, so this algorithm may not be applicable in our situation.

Finally, we note that it seems natural to use (for hashing purposes) pairs of matrices with smaller entries, like $(A(1), B(1))$ or $(A(1), B(-1))$. However, the matrices $(A(1), B(-1))$ do not generate a free semigroup (they satisfy the ``braid relation" $ABA=BAB$), so our approach to collision resistance does not work in this case. The matrices $(A(1), B(1))$ do generate a free semigroup, but they generate the {\it whole} semigroup of matrices in $SL_2(\Z)$ with nonnegative entries. This property was used in \cite{TZ_attack}  to arrange an attack that was successful in recovering a preimage of a hash, thereby showing that the corresponding Cayley hash function is not {\it preimage resistant}.

We leave it at that.

{\renewcommand{\baselinestretch}{1.17}


\begin{thebibliography}{ABC}


\bibitem{AB}
G. Arzhantseva and A. Biswas, {\it  Logarithmic girth expander graphs of $SL_n(\F_p)$}, J. Algebraic Combinatorics {\bf 56} (2022), 691--723.

\bibitem{AO}
G.~Arzhantseva and A.~Olshanskii, {\it Genericity of the class of groups in which subgroups with a lesser number of generators are free,} (Russian) Mat. Zametki \textbf{59} (1996), no. 4, 489--496.

%\bibitem{Baumslag}
%B. Baumslag and S. J. Pride, {\it Groups with two more generators than relators}, J. London Math. Soc. (2)  {\bf 17} (1978), 425-–426.

\bibitem{bookof10}
F. Bassino, I. Kapovich, M. Lohrey, A. Myasnikov, C. Nicaud, A. Nikolaev, I. Rivin, V. Shpilrain, A. Ushakov, P. Weil, {\it  Complexity and Randomness in Group Theory: GAGTA Book 1}. Walter de Gruyter, 2020.

\bibitem{Problems}
G. Baumslag, A. G. Myasnikov, V. Shpilrain, {\it Open problems in
combinatorial group theory},\\
\url{http://shpilrain.ccny.cuny.edu/gworld/problems/oproblems.html}

\bibitem{Birget}
J.-C. Birget, A. Yu. Olshanskii, E. Rips, M. V. Sapir, {\it Isoperimetric functions of groups and computational complexity of the word problem}, Ann. of Math. (2) {\bf 156} (2002), 467--518.

\bibitem{BKL}
J. Birman, K. H. Ko, S. J. Lee, {\it A new approach to the word and conjugacy problems in the braid groups}, Adv. Math. {\bf 139} (1998), 322--353.

\bibitem{BG}
J. Bourgain, A. Gamburd, {\it Uniform expansion bounds for Cayley
graphs of $SL_2({\F}_p)$}. Ann. of Math. (2) {\bf 167}
(2008),  625--642.

\bibitem{BSV}
L. Bromberg, V. Shpilrain, A. Vdovina, {\it
Navigating in the Cayley graph of $SL_2(\F_p)$ and applications to
hashing}, Semigroup Forum {\bf 94} (2017), 314--324.

\bibitem{Carnevale}
A. Carnevale, M. Cavaleri, {\it Partial word and equality problems and Banach
densities}, Adv. Math. {\bf  368} (2020), 107133.

\bibitem{languages}
T. Ceccherini-Silberstein, W. Woess, {\it  Growth and ergodicity of context-free languages}, Trans. Amer. Math. Soc. {\bf  354} (2002), 4597--4625.

\bibitem{Cohen}
J. M. Cohen, {\it Cogrowth and amenability of discrete groups},  J. Funct. Anal. {\bf 48} (1982), 301--309.

\bibitem{Grigorchuk}
R. I. Grigorchuk, {\it Symmetrical random walks on discrete groups},
in: Advances in Probability and Related Topics, Vol. {\bf 6}, pp. 285--325, Marcel Dekker 1980.

%\bibitem{graphs}
%D. Grigoriev and  V. Shpilrain, {\it Authentication schemes from actions on graphs, groups, or rings}, Ann. Pure Appl. Logic. {\bf 162} (2010), 194--200.

%\bibitem{tropical}
%D. Grigoriev and  V. Shpilrain, {\it Tropical cryptography},  Comm. Algebra {\bf 42} (2014), 2624--2632.

\bibitem{Gromov}
M. Gromov, {\it Asymptotic invariants of infinite groups}, Geometric group theory, Vol. 2 (Sussex, 1991), 1--295. London Math. Soc. Lecture Note Ser., {\bf 182}
Cambridge University Press, Cambridge, 1993.

\bibitem{Guba}
V. Guba, {\it Amenability problem for Thompson's group $F$: state of the art}, Groups, Complexity, Cryptology {\bf 15} (2023.

%\bibitem{secret}
%M. Habeeb,  D. Kahrobaei, V. Shpilrain, {\it A secret sharing scheme based on group presentations and the word problem}, Contemp. Math., Amer. Math. Soc. {\bf 582} (2012),
%143--150.

%\bibitem{semidirect}
%M. Habeeb, D. Kahrobaei, C. Koupparis, V. Shpilrain, {\it Public key
%exchange using semidirect product of (semi)groups},  in: ACNS 2013, Lecture Notes Comp. Sc. {\bf 7954} (2013), 475--486.

\bibitem{Harvey}
D. Harvey and J. van der Hoeven, {\it Integer multiplication in time $O(n \log n)$}, Ann. of Math. {\bf 193} (2021), 563--617.

\bibitem{H}
H. A. Helfgott,{\it Growth and generation in $SL_2(\Z/p\Z$)} Ann. of
Math. (2) {\bf 167} (2008),  601--623.

\bibitem{RamaNotices}
D. Kahrobaei, R. Flores, M. Noce, {\it Group-based cryptography in the quantum era}, Notices Amer. Math. Soc. {\bf 70} (2023), 752--763.

\bibitem{bookDK}
D. Kahrobaei, R. Flores, M. Noce, M. Habeeb,  C. Battarbee,  {\it Applications of Group Theory in Cryptography: Post-quantum Group-based Cryptography}, Amer. Math. Soc. Surveys and Monographs, 2023, to appear.

\bibitem{KMSS}
I. Kapovich, A. G. Myasnikov, P. Schupp, V. Shpilrain,
{\it  Generic-case complexity, decision problems in group theory and random
walks}, J.  Algebra {\bf 264} (2003), 665--694.

\bibitem{KMSS2}
I. Kapovich, A. G. Myasnikov, P. Schupp, V. Shpilrain, {\it
Average-case complexity and decision problems in group theory},
 Advances in Math. {\bf 190} (2005), 343--359.

%\bibitem{KRSS}
%I. Kapovich, I. Rivin,  P. Schupp, V. Shpilrain, {\it Densities in
%free groups and $\mathbb{Z}^k$, visible points and test elements},
%Math. Res. Lett. {\bf 14}  (2007),  263--284.

%\bibitem{KSS}
%I. Kapovich, P. Schupp, and V. Shpilrain, {\it Generic properties of
%Whitehead's algorithm and isomorphism rigidity of random one-relator
%groups}, Pacific J. Math. {\bf 223} (2006), 113--140.

\bibitem{OK}
O. G.  Kharlampovich, {\it A finitely presented solvable group with unsolvable word problem.} (Russian) Izv. Akad. Nauk SSSR Ser. Mat. {\bf  45} (1981), 852--873.

\bibitem{KMS}
O. G.  Kharlampovich, A. G. Myasnikov, M. Sapir, {\it Algorithmically complex residually finite groups}, Bull. Math. Sci. {\bf 7} (2017),  309--352.

\bibitem{Knuth}
D. E. Knuth, {\it The analysis of algorithms}. Actes du Congr\`es International des Math\'ematiciens (Nice, 1970), Tome 3, pp. 269–-274. Gauthier-Villars, Paris, 1971.

\bibitem{Levin}
L. Levin, {\it Average case complete problems}, SIAM J. Comput. {\bf 15} (1986),  285--286.

\bibitem{Lohrey}
M. Lohrey, {\sl The compressed word problem for groups}, Springer Briefs Math., Springer, New York, 2014.

\bibitem{LS}
R.~Lyndon and P.~Schupp, \emph{Combinatorial Group Theory,} Springer-Verlag, 1977. Reprinted in the ``Classics in  mathematics'' series, 2000.

\bibitem{Mattes}
C. Mattes, A. Weiss, {\it Improved Parallel Algorithms for Generalized Baumslag
Groups}, in: LATIN 2022: Theoretical Informatics, Lecture Notes in Comput. Sci. {\bf 13568} (2022), 658--675.


\bibitem{Vershik}
A. G. Myasnikov, V. Roman’kov, A. Ushakov, and A. Vershik, {\it The word and geodesic problems in free solvable groups}, Trans. Amer. Math. Soc. {\bf 362} (2010), 4655--4682.
%\bibitem{MSU2}
%A. G. Myasnikov, V. Shpilrain, and A. Ushakov, {\sl Non-commutative
%cryptography and complexity of group-theoretic problems}, Amer.
%Math. Soc. Surveys and Monographs, 2011.


\bibitem{book2}
A. G. Myasnikov, V. Shpilrain, and A. Ushakov, {\it Group-based cryptography}. Birkh\"auser Verlag, Basel-Boston-Berlin, 2008.

\bibitem{book3}
A. G. Myasnikov, V. Shpilrain, and A. Ushakov, {\it Non-commutative cryptography and complexity of group-theoretic problems}. Amer. Math. Soc. Surveys and Monographs, 2011.

\bibitem{MU}
A. G. Miasnikov, A. Ushakov, {\it Random van Kampen diagrams and algorithmic
problems in groups},  Groups Complex. Cryptol. {\bf 3} (2011),
121--185.

%\bibitem{Nikolaev}
%I. Nikolaev, {\it Mapping class groups are linear}, in: Topological algebras and their applications, 193--200,
%De Gruyter Proc. Math., De Gruyter, Berlin, 2018.

\bibitem{OSh}
A. Olshanskii, V. Shpilrain, {\it Linear average-case complexity of algorithmic problems in groups,} preprint. %{\url https://arxiv.org/abs/2205.05232}


\bibitem{PatRaz}
M. S. Paterson, A. A. Razborov, {\it The set of minimal braids is co-NP-complete}, J. Algorithms {\bf 12} (1991), 393--408.

\bibitem{Pe}
D.~Peifer, \emph{Artin groups of extra-large type are biautomatic,} J. Pure Appl. Algebra  \textbf{110}  (1996), 15--56.

\bibitem{Pollicott}
M. Pollicott, {\it Maximal Lyapunov exponents for random matrix products},
Invent. Math. {\bf 181} (2010), 209--226.

%\bibitem{Roig}
%A. Roig, E. Ventura, P. Weil,  {\it On the complexity of the Whitehead
%minimization problem}, Internat. J. Algebra Comput. {\bf  17} (2007),  1611--1634.

\bibitem{Roy}
M. Roy, E. Ventura, P. Weil, {\it The central tree property and algorithmic problems
on subgroups of free groups}, J. Group Theory, to appear.

%\bibitem{Sharp}
%R. Sharp, {\it Local limit theorems for free groups}, Math. Ann. {\bf 321} (2001), 889--904.

\bibitem{Saul}
S. Schleimer, {\it Polynomial-time word problems}, Comment. Math. Helv. {\bf 83} (2008), 741--765.

%\bibitem{witness}
%V. Shpilrain, {\it Search and witness  problems
%in group theory}, Groups, Complexity, Cryptology {\bf 2} (2010),
%231--246.

\bibitem{Sublin}
V. Shpilrain, {\it Sublinear time algorithms in the theory of groups
and semigroups}, Illinois J. Math. {\bf 54} (2011), 187--197.

\bibitem{White}
V. Shpilrain, {\it  Average-case complexity of the Whitehead problem for free groups}, Comm. Algebra {\bf 51} (2023), 799--806. %{\url https://arxiv.org/abs/2105.01366}

%\bibitem{SU}
%V. Shpilrain and A. Ushakov, {\it Thompson's group and public key cryptography}, in  ACNS 2005, Lecture Notes Comp. Sc. {\bf 3531} (2005), 151--164.
%
%\bibitem{Ushakov}
%A. Ushakov, {\it Algorithmic theory of free solvable groups: randomized computations},  J. Algebra {\bf 407} (2014), 178--200.

\bibitem{Sisto}
A. Sisto, {\it Tracking rates of random walks},  Israel J. Math. {\bf 220}
(2017), 1--28.

\bibitem{TZ_attack}
J.-P. Tillich and G. Z\'emor,  {\it Group-theoretic hash functions},
in Proceedings of the First French-Israeli Workshop on Algebraic
Coding, Lecture notes  Comp. Sci. {\bf  781} (1993),   90--110.

\bibitem{Timoshenko}
E. I. Timoshenko, {\it Certain algorithmic questions for metabelian groups}, Algebra and Logic {\bf 12} (1973), 132--137.


\bibitem{Wang95}
J.~Wang, \emph{Average-case completeness of a word problem in
groups,} Proc. of the 27-th Annual Symposium on Theory of
Computing, ACM Press, New York, 1995, 325--334.

\bibitem{Wh}
J.~H.~C.~Whitehead, {\it On equivalent sets of elements in free
groups}, Ann. of Math. \textbf{37} (1936), 782--800.

\bibitem{Woess}
W. Woess, {\sl Random walks on infinite graphs and groups}, Cambridge Tracts in Mathematics, {\bf 138}. Cambridge University Press, Cambridge, 2000.

\bibitem{Young}
R. Young, {\it Averaged Dehn functions for nilpotent groups}, Topology
{\bf 47} (2008),  351--367.


\end{thebibliography}
\end{document}